\documentclass[preprint,12pt]{elsarticle}

\usepackage{amssymb}
\usepackage{amsfonts}
\usepackage{amsmath}
\usepackage[T1]{fontenc}
\usepackage{graphicx}
\usepackage[mathscr]{eucal}
\usepackage[absolute]{textpos}
\usepackage{amsmath, amscd}

\journal{Topology and its Applications}

\begin{document}

\begin{frontmatter}


\title{Milnor-Thurston homology groups of the Warsaw Circle}
\author{Janusz Przewocki}
\ead{j.przewocki@gmail.com}
\ead[url]{http://www.impan.pl/~jprzew/}
\address{Institute of Mathematics of the Polish Academy of Sciences, 8 Sniadeckich st., 00-956 Warsaw, Poland}





\newtheorem{df}{Definition}
\newtheorem{thm}[df]{Theorem}
\newtheorem{lem}[df]{Lemma}
\newtheorem{cor}[df]{Corollary}

\begin{abstract}
Milnor-Thurston homology theory is a construction of homology theory that is based on measures. It is known to be equivalent to singular homology theory in case of manifolds and complexes. Its behaviour for non-tame spaces is still unknown. This paper provides results in this direction. We prove that Milnor-Thurston homology groups for the Warsaw Circle are trivial except for the zeroth homology group which is uncountable dimensional. Additionally, we prove that the zeroth homology group is non-Hausdorff for this space with respect a natural topology that was proposed by Berlanga.

\end{abstract}


\end{frontmatter}


\section{Introduction}

Homology theory is one of the main tools used in algebraic topology. We can find its origins at the end of 19th century when people started their research in order to understand manifolds and, more generally, simplicial complexes. Later, by creating singular homology theory, the concepts had been generalized to all topological spaces. Singular homology theory is well understood for tame spaces (e.g. simplicial complexes, manifolds) and some steps to grasp its behaviour for more complicated spaces has been taken by a new emerging field, which might be described as ``algebraic topology of non-tame spaces'' (as it was called in \cite{AZNontame} or, more briefly, as ``wild algebraic topology'', as it was called in \cite{KEdaTalk}; this field is focused on spaces like: Hawaiian Earring \cite{MorgMor} \cite{CanCon} \cite{KEda}, Harmonic Archipelago \cite{BogSier} \cite{DRepovs}, Sierpinski-gasket \cite{SierpGask}, Griffiths' space \cite{Grif}, etc.).

The previously mentioned Hawaiian Earring whose first singular homology group has been algebraically described in \cite{KEda} is a sequence of planar circles tangential to each other at the single point which have its diameters converging to zero. The reason why this fairly simple space has a very complicated homology group is that chains in singular homology consist of a finite number of simplices, and we clearly see that the Hawaiian Earring has infinite number of ``cycles'' (and unlike in the case of a one point union of countably many circles with CW-topology there exist paths that wind around infinitely many circles). 
 
As a consequence, a homology theory with chains allowing an infinite numbers of simplices may be more suitable for the Hawaiian Earring or, more generally, for wild topological spaces. The aim of this paper is to study such a theory: \emph{Milnor-Thurston homology theory}.  

Milnor-Thurston homology theory (which is also called \emph{measure homology theory}) is a construction where chains are Borel measures on the space of simplices (equipped with the compact-open topology). This homology theory was introduced by Thurston in \cite[Chapter 6.]{Thurst}. The reason why we should use measures instead of finite linear combinations of simplices is that it gives us more freedom to perform certain constructions. A good example we can observe in Thurston's proof of Gromov's theorem which roughly says that the hyperbolic volume of a hyperbolic manifold is a topological invariant (see proof of Theorem 6.2 in \cite{Thurst}). Also the proof of Mostow Ridgity Theorem as presented in Ratcliffe book \cite{Ratcl} uses this homology theory in a smooth setting (this type of homology was proved by L\"oh to be isometrically equivalent to the one we use here \cite{Claraloh}).  
 
A formal description and some basic properties of Milnor-Thurston homology theory can be found in \cite{AZastrow} or \cite{SKHansen}. Additionally, these homology groups can be equipped with a natural topological structure that we call here ``Berlanga topology'' (see \cite{Berla}). 

One of the basic properties of Milnor-Thurston homology theory is that it coincides with singular homology theory for CW-complexes \cite{AZastrow}. Currently there are no results on behaviour of the theory for non-tame spaces. Therefore, we decided to investigate the case of one of the basic wild spaces -- the Warsaw Circle. 

In Section \ref{sec:basics} we recall basic facts on Milnor-Thurston homology theory. In Section \ref{sec:groups} we prove that nonzero order Milnor-Thurston homology groups for the Warsaw Circle are trivial, and that the zeroth Milnor-Thurston homology group is an uncountable dimensional vector space.

The reason why we took the Warsaw Circle is that it is one of the simplest example of a non-triangulable space. Additionally, one can expect that the results should be different than in the singular homology case, since our theory admits infinite cycles. In particular, the reasonable question may be: ``Does the Milnor-Thurston homology theory detect the circular shape of the Warsaw Circle?''. Being more precise, we ask whether $\mathcal{H}_1(W)$ equals $\mathbb{R}$. As we shall see, it is not the case.

Finally, in Section \ref{sec:berlanga1} we define Berlanga topology and we prove that it is non-Hausdorff for the zeroth Milnor-Thurston homology group of the Warsaw Circle. This answers in a negative way the question posed in \cite{Berla}.

\section{Basics of Milnor-Thurston homology theory}
\label{sec:basics}

Milnor-Thurston homology theory is a version of the homology theory which admits chains with infinite number of simplices. As was mentioned before, it has first appeared in connection with hyperbolic geometry in late 70s. Two decades later it was formalised by Zastrow \cite{AZastrow} and Hansen \cite{SKHansen} independently, and its definition was generalized to all topological spaces. 

It has been proved that this theory satisfies the Eilenberg-Steenrod axioms with weakened Excision Axiom, which is equivalent to the classical one for well behaved spaces (including all normal spaces) \cite{AZastrow}. This fact is of a great importance. First of all, it implies that for tame spaces the theory is the same as singular homology theory. Additionally, it yields the Mayer-Vietoris theorem. 

In this paper we use calligraphic letters ($\mathcal{C}$, $\mathcal{H}$, etc.) for constructions in Milnor-Thurston homology theory and ordinary letters for the corresponding constructions in singular homology theory ($C$, $H$, etc.). Accordingly, the Mayer-Vietoris theorem for Milnor-Thurston homology theory asserts
\begin{thm}
Let $X$ be a normal topological space and let $A$ and $B$ be open subspaces such that $X = A \cup B$. Then for all $n \geq 0$ there exists a natural homomorphism $\partial_*: \mathcal{H}_{n}(X) \rightarrow \mathcal{H}_{n-1}(A \cap B)$ such that the following sequence is exact:
\begin{multline}
\begin{CD}
\cdots @>\text{$(i_{* n}, j_{* n})$}>>  \mathcal{H}_{n}(A) \oplus \mathcal{H}_{n}(B) @>\text{$s_{* n} - t_{* n}$}>> \mathcal{H}_{n}(X)  @>\text{$\partial_*$}>> \mathcal{H}_{n-1}(A \cap B) @>>>
\end{CD}  \\
\begin{CD}
\cdots \rightarrow  \mathcal{H}_{0}(A \cap B) @>\text{$(i_{* 0}, j_{* 0})$}>> \mathcal{H}_{0}(A) \oplus \mathcal{H}_{0}(B)  @>\text{$s_{* 0} - t_{* 0}$}>> \mathcal{H}_{0}(X) @>>> 0
\end{CD}  \nonumber
\label{eq:diagram}
\end{multline}
where $i: A \cap B \rightarrow A$,  $j: A \cap B \rightarrow B$, $s: A \rightarrow X$, $t: B \rightarrow X$ are inclusion maps. 
\label{thm:mayervietoris}
\end{thm}

Let $X$ be a topological space. We shall start construction of Milnor-Thurston homology theory with defining its chain complex $\mathcal{C}_*(X)$. 

A \emph{signed measure} is a $\sigma$-additive set function with possibly negative values such that it is zero on the empty set. Every measure considered here is a signed measure therefore we shall simply call them  \emph{measures}.  

We will be concerned with singular simplices (continuous functions from the standard simplex $\Delta^k$ to $X$, where $k$ is a non-negative integer) as in the case of the singular homology theory. We endow the space of singular simplices $C^0(\Delta^k, X)$ with the \emph{compact-open topology}. On this space we consider Borel sets $\mathbb{B}(C^0(\Delta^k, X))$ -- the smallest $\sigma$-algebra generated by open sets. A measure defined for all Borel sets on the given space is called a \emph{Borel measure}. It is said to be finite when it has finite values for all Borel sets. A \emph{carrier} of Borel measure $\mu$ is a set $D \subset C^0(\Delta^k, X) $ such that all measurable subsets of $C^0(\Delta^k, X) \setminus D$ are zero sets (i.e. sets such that every Borel subset has measure zero). Notice, that for a given measure carrier is not a uniquely determined set. 

Now, we define the real vector space $\mathcal{C}_k(X)$ as the set of finite Borel measures on $C^0(\Delta^k, X)$ with some compact carrier. Next, we will construct boundary homomorphisms in order to endow the sequence $\mathcal{C}_k(X)$ with a structure of chain complex.

Given a continuous function $f: C^0(\Delta^k, X) \rightarrow C^0(\Delta^l, Y)$ where $Y$ is a topological space and $l$ is a non-negative integer and an arbitrary Borel measure $\mu$ on $C^0(\Delta^k, X)$  we define the \emph{image measure} $f\mu$ on $C^0(\Delta^l, Y)$ with the formula:
\begin{displaymath}
(f \mu) (A) = \mu(f^{-1}(A)), \quad \textrm{ for } A \in \mathbb{B}(C^0(\Delta^l, Y)).
\end{displaymath} 
Now, we construct the boundary operator $\partial$ in the usual way:
\begin{displaymath}
\partial = \sum_{i = 0}^k (-1)^i \partial_i, 
\end{displaymath}
where $\partial_i$ sends a measure to the image measure under the map $\sigma \mapsto \sigma \circ \delta_i$ with $\delta_i$ being the usual inclusion of $\Delta^{k-1}$ as a face of $\Delta^k$. We can prove \cite[Corollary 2.9]{AZastrow} that $\mathcal{C}_*(X)$ with the boundary operator is a chain complex. 

The Milnor-Thurston homology groups $\mathcal{H}_*(X)$ are then defined as homology groups of this chain complex $\mathcal{C}_*(X)$. Additionally, $\mathcal{C}_*$ can be treated as a functor from the category of topological spaces to the category of chain complexes. Thus, we can define relative homology groups $\mathcal{H}_*(X, A)$ in a natural way.

\section{Milnor-Thurston homology groups for the Warsaw Circle}
\label{sec:groups}

\begin{figure}
\begin{center}
\includegraphics[height=0.25\textheight]{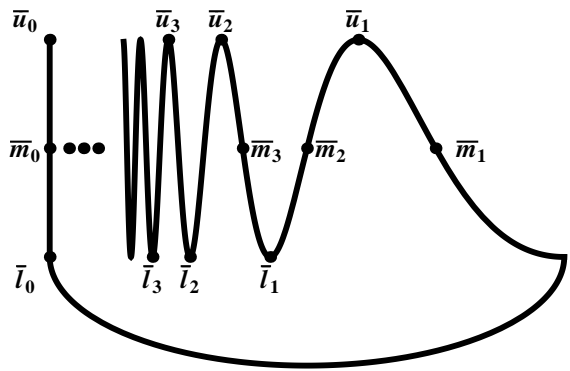}

\caption{The Warsaw Circle with distinguished points}
\label{fig:WarsawCircle}
\end{center}
\end{figure}

The Warsaw Circle (see Figure \ref{fig:WarsawCircle}) is a well known space that serves as a counterexample in many cases (e.g. it is a basic example of a non-triangulable space). We define it as a subset of $\mathbb{R}^2$ that consists of:
\begin{itemize}
    \item  the part of Topologist Sine Curve $\{(x, y) \in \mathbb{R}^2 \mid y = \sin 1/x \}$ between the line $x = 0$ and the rightmost minimum, 
    \item the ``accumulation line'' $\{(0, y) \in \mathbb{R}^2 \mid -1 \leq y \leq 1  \}$, 
    \item an arc connecting the point $(0, -1)$ with the rightmost minimum.
\end{itemize}

Throughout this paper $W$ will denote the Warsaw Circle. We distinguish three families of points $\{\overline{l}_k\}_{k=0}^\infty$, $\{\overline{m}_k \}_{k=0}^\infty$, $\{\overline{u}_k\}_{k=0}^\infty$ in $W$ (see Figure \ref{fig:WarsawCircle}). Being more precise, $\overline{u}_1$ is a first maximum of the sinusoid to the left of the rightmost minimum. Next, $\overline{u}_2$ is a first maximum of the sinusoid to the left of $\overline{u}_1$ and so on. Finally, $\overline{u}_0 = (0, 1) \in \mathbb{R}^2$. The families  $\{\overline{m}_k \}_{k=0}^\infty$, $\{\overline{l}_k\}_{k=0}^\infty$ are defined in an analogous manner. 
 
In order to compute Milnor-Thurston homology groups for $W$, we shall cover it by two open subsets $L$ and $U$. Let $L$ be an intersection of $W$ with the half-plane $\{(x, y) \in \mathbb{R}^2 \mid y < \eta\}$, where $0 < \eta < 1$. Similarly, let $U$ be an intersection of $W$ with  $\{(x, y) \in \mathbb{R}^2 \mid y > -\eta\}$.

Since Milnor-Thurston homology groups in this paper are conveniently described by infinite sequences of numbers we shall use a special notation. So let us assign some element $x_k$ to any non-negative integer $k$. The \emph{sequence} of such elements shall be denoted by $x_\bullet$ (usually we will be concerned with sequences of real numbers). 

Now, we will use the Mayer-Vietoris theorem to do the calculations, but first, in order to deal with the variety of possible singular simplices, we shall need an auxiliary result.  So let $S$ denote a convergent sequence $s_\bullet$ with its limit $s_0$. In this context observe
\begin{equation}
U \simeq S
\label{eq:homotop1} 
\end{equation}
\begin{equation}
L \simeq S 
\label{eq:homotop2}
\end{equation}
 
Since homology groups are homotopy invariant the next lemma allows us to calculate them for $U$, $L$ and $U \cap L$. 
\begin{lem}
If $n > 0$, then $\mathcal{H}_n(S) = 0$ and $\mathcal{H}_0(S) \cong \ell^1$. Here $\ell^1$ denotes the vector space of absolutely summable sequences. 
\label{lem:sequence}
\end{lem}

\textbf{Remark.} Note that Milnor-Thurston homology groups are simply real vector spaces, so the lemma states that $\mathcal{H}_0(S) \cong \bigoplus_{\mathfrak{c}} \mathbb{R}$. Nevertheless, we shall use more concrete description of that vector space, namely $\ell^1$ (we adapted that notation from functional analysis, however we do not treat it as a Banach space). The reason is that the Milnor-Thurston homology classes that appear in this paper usually have a simple description as measures concentrated on distinguished points and coefficients of these measures form absolutely summable infinite sequences.  

\textbf{Proof.} We can see that
\begin{displaymath}
C^0(\Delta^n, S) = \{s_k^n: \Delta^n \rightarrow S \mid s_k^n \textrm{ sends any point of $\Delta^n$ to $s_k$}, k \geq 0 \}.
\end{displaymath} 

The topology of $C^0(\Delta^n, S)$ is the same as for $S$ because $\{s_k^n\}_{k=0}^\infty$ is a convergent sequence with limit $s_0^n$. In addition, every set of this space is Borel and every two Borel measures that are equal on singletons $\{s_k^n \}$ for every $k$ are equal.  Therefore, we can identify a sequence of real numbers $a_\bullet$ with a measure $\mu_a$ in the way that $\mu_a(\{s_k^n \}) = a_k$ for every non-negative integer $k$. Additionally, for every sequence $a_\bullet$ the measure $\mu_a$ has a compact carrier (that is the whole space $C^0(\Delta^n, S)$) and $\mu_a$ is finite if and only if 
\begin{displaymath}
\sum_{k = 0}^\infty |a_k| < \infty.
\end{displaymath}
Consequently,
\begin{displaymath}
\mathcal{C}_n (S) \cong \ell^1 := \{a_\bullet \mid a_k \in \mathbb{R},  \sum_{k = 0}^\infty |a_k| < \infty \}.
\end{displaymath} 
We shall study the boundary operator: $\partial: \mathcal{C}_n(S) \rightarrow \mathcal{C}_{n-1}(S)$ for $n > 0$. Both $\mathcal{C}_n(S)$ and $\mathcal{C}_{n-1}(S)$ shall be identified with $\ell^1$. We have that $\partial_i a_\bullet = a_\bullet$ (recall that the definition of $\partial_i$ was given in Section \ref{sec:basics}). To see this reader should notice that $\partial_i s_k^n = s_k^{n-1}$. As a consequence,
\begin{displaymath}
\partial a_\bullet = \sum_{i = 0}^n (-1)^i \partial_i a_\bullet = a_\bullet \cdot \sum_{i = 0}^n (-1)^i.  
\end{displaymath}  
From here, $\partial = 0$ when $n$ is odd and $\partial = \mathrm{id}$ when $n$ is even. Thus, homology spaces are trivial for $n > 0$. 

On the other hand, we have $\partial = 0$, for $n = 0$. Hence, every element in $\mathcal{C}_0(S)$ is a cycle. Because  $\partial = 0$, for $n = 1$, there are no boundaries and $\mathcal{H}_0(S) = \mathcal{C}_0(S) \cong \ell^1$. 
\begin{flushright} $\square$ \end{flushright} 

Finally, using the Mayer-Vietoris theorem we can calculate the homology groups. 
\begin{thm}
If $n > 0$, then $\mathcal{H}_n(W) = 0$. 
\label{thm:main1}
\end{thm}

\textbf{Proof.} The Mayer-Vietoris sequence
\begin{multline}
\begin{CD}
\cdots @>\text{$(i_{* n}, j_{* n})$}>>  \mathcal{H}_{n}(U) \oplus \mathcal{H}_{n}(L) @>\text{$s_{* n} - t_{* n}$}>> \mathcal{H}_{n}(W)  @>\text{$\partial_*$}>> \mathcal{H}_{n-1}(U \cap L) @>>>
\end{CD}  \\
\begin{CD}
\cdots \rightarrow  \mathcal{H}_{0}(U \cap L) @>\text{$(i_{* 0}, j_{* 0})$}>>\mathcal{H}_{0}(U) \oplus \mathcal{H}_{0}(L)  @>\text{$s_{* 0} - t_{* 0}$}>> \mathcal{H}_{0}(W) @>>> 0
\end{CD}  \nonumber
\label{eq:diagram}
\end{multline}
is exact (cf. Theorem \ref{thm:mayervietoris}). 
Then by equations (\ref{eq:homotop1}), (\ref{eq:homotop2}) and Lemma \ref{lem:sequence} we see that the Mayer-Vietoris sequence has zeros on the left hand side of $\mathcal{H}_1(W)$. Thus, higher homology groups are trivial and the only case to investigate is where $n = 1$.  Again, by exactness of the Mayer-Vietoris sequence we see that $\partial_*: \mathcal{H}_1(W) \rightarrow \mathcal{H}_0(U \cap L)$ is a monomorphism. As a result, 
\begin{displaymath}
\mathcal{H}_1(W) \cong \mathrm{ker}(i_{* 0}, j_{* 0}). 
\end{displaymath}
Therefore, we need to find the kernel of $(i_{* 0}, j_{* 0})$.

Again, by equations (\ref{eq:homotop1}), (\ref{eq:homotop2}) 
and Lemma \ref{lem:sequence} we can see that
\begin{displaymath}
\mathcal{H}_0(U) \cong \mathcal{H}_0(L) \cong  \mathcal{H}_0(U\cap L) \cong  \mathcal{C}_0(S) \cong \ell^1,
\end{displaymath}
so from now on we shall identify elements of all these homology groups with absolutely summable real sequences.

So, let $m_\bullet \in \ell^1$ denote some homology class of $\mathcal{H}_0(U \cap L)$ that is represented by a measure supported on respective points $\overline{m}_k$ (the numbers $m_k$ are the values of the measure on the singletons). Similarly, a homology class in $\mathcal{H}_0(U)$ shall be identified with some $u_\bullet \in \ell^1$ and a homology class in $\mathcal{H}_0(L)$ shall be identified with some $l_\bullet \in \ell^1$. This will allow us to write down equations for $i_{* 0}$ and $j_{* 0}$.

In order to investigate $i_{* 0}$, we have to associate a measure supported on $\overline{u}_\bullet$ with a measure supported on $\overline{m}_\bullet$ that represents the same homology class in $U$. So, let $\mu$ be a measure supported on $\overline{m}_\bullet$ (cf. Figure \ref{fig:WarsawCircle}) represented by the sequence $m_\bullet$. We will construct a measure supported on $\overline{u}_\bullet$ which belongs to the same $\mathcal{H}_0(U)$-homology class as $\mu$. Let $\sigma_0$ be a singular 1-simplex that connects  $\overline{m}_0$ with $\overline{u}_0$. Next, let $\sigma_{2k}$ be a singular 1-simplex that connects $\overline{m}_{2k}$ with $\overline{u}_k$ and let $\sigma_{2k+1}$ be a singular 1-simplex that connect $\overline{m}_{2k+1}$ with $\overline{u}_k$. Now, let $\nu = \sum_{k=0}^\infty m_k \delta_{\sigma_k}$, where $\delta_{\sigma_k}$ is the Dirac measure supported on $\sigma_k$. We can see, that $\nu \in \mathcal{C}_1(U)$, since $\nu$ is finite and has a compact carrier (because $\sigma_\bullet$ is a convergent sequence). The measure $\mu + \partial \nu$ is supported on $\overline{u}_\bullet$, its coefficients depend on $m_\bullet$ as described below. From the definition of $\sigma_0$ we have that          
\begin{equation}
u_0 = m_0. 
\label{eq:zero}
\end{equation}
Furthermore, from the definitions of $\sigma_{2k}$ and $\sigma_{2k+1}$ we have
\begin{equation}
u_k = m_{2k} + m_{2k-1} \quad \textrm{for } k > 0.
\label{eq:one}
\end{equation}
These are equations that describe $i_{* 0}$. 

In the similar way, we can write down equations for $j_{* 0}$: 
\begin{equation}
l_k = m_{2k} + m_{2k+1}.
\label{eq:two}
\end{equation}

We will describe $(i_{* 0}, j_{* 0})$ in more elegant way. So let $x_{2k} = u_k$ and $x_{2k-1} = l_k$. From now on an absolutely summable sequence $x_\bullet$ shall be identified with elements of $\mathcal{H}_0(U) \oplus \mathcal{H}_0(L)$. In this notation, equations (\ref{eq:zero}), (\ref{eq:one}), (\ref{eq:two}) yield
\begin{equation}
x_k = \left\{ \begin{array}{lr}
m_0 & \textrm{for } k = 0, \\
m_k + m_{k-1} & \textrm{for } k > 0. 
\end{array} \right.
\label{eq:basic}
\end{equation}

Now, we see that the kernel of $(i_{* 0}, j_{* 0})$ and, consequently, $\mathcal{H}_1(W)$ are trivial. 
\begin{flushright} $\square$ \end{flushright}

\begin{thm}
The vector space $\mathcal{H}_0(W)$ is uncountable-dimensional.
\label{thm:undim}
\end{thm}

\textbf{Proof.} We shall use results from the proof of Theorem \ref{thm:main1}; mostly equation (\ref{eq:basic}).  Once again we will use the Mayer-Vietoris theorem to do the calculations. The Mayer-Vietoris sequence is (cf. Theorem \ref{thm:mayervietoris})
\begin{multline}
\begin{CD}
0 \rightarrow  \mathcal{H}_{0}(U \cap L) @>\text{$(i_{* 0}, j_{* 0})$}>>\mathcal{H}_{0}(U) \oplus \mathcal{H}_{0}(L)  @>\text{$s_{* 0} - t_{* 0}$}>> \mathcal{H}_{0}(W) @>>> 0. 
\end{CD}  \nonumber
\label{eq:diagram}
\end{multline}

Now we can see that $\mathcal{H}_0(W)$ is a quotient $\ell^1 / h(\ell^1)$ where $h: \ell^1 \rightarrow \ell^1$ is the map defined by equation ($\ref{eq:basic}$). This equation can be inverted so that, given an arbitrary sequence $x_\bullet$, we can find unique numbers $m^x_k$ that satisfy it:
\begin{equation}
m_k^x = \sum_{i = 0}^k (-1)^{i+k} x_i.
\label{eq:solution}
\end{equation} 
An element $x_\bullet \in \ell^1$ represents a nonzero homology class in $\mathcal{H}_0(W)$ if it is not in the image of $(i_{* 0}, j_{* 0})$ or, equivalently, when the corresponding $m^x_\bullet$ is not an absolutely summable sequence. 

For any space $X$ there is the natural inclusion of singular chains into Milnor-Thurston chains: $C_k(X; \mathbb{R}) \rightarrow \mathcal{C}_k(X)$ \cite[p. 389]{AZastrow}. It induces a homomorphism on the level of homology. Thus, we can form the following definition:
\begin{df}
A homology class in $\mathcal{H}_k(X)$ shall be called \emph{singular homology class} if it lies in the image of $H_k(X; \mathbb{R}) \rightarrow \mathcal{H}_k(X)$. Otherwise it shall be called \emph{non-singular homology class}. 
\end{df}

Next, we shall find a one dimensional subspace of $\mathcal{H}_0(W)$ which corresponds to singular homology classes. In singular homology theory we consider chains with only finite numbers of simplices, so we will restrict ourselves to considering sequence $x_\bullet$ with finitely many nonzero elements. We will prove that such an element $x_\bullet \in \ell^1$ represents the same homology class as $y_\bullet \in \ell^1$ of the form $y_\bullet = (\alpha, 0, 0, 0, \dots)$, for some $\alpha \in \mathbb{R}$. Let $N$ denote the biggest index of nonzero element in $x_\bullet$ then for $k > N$ we have
\begin{displaymath}
m_k^{x-y} = (-1)^k \left(\sum_{i = 0}^N (-1)^i x_i - \alpha\right).
\end{displaymath}
So putting $\alpha = \sum_{i = 0}^N (-1)^i x_i$, yields $m_k = 0$. Thus, it is absolutely summable and $x_\bullet - y_\bullet$ represents the zero homology class. 

Now, we shall prove that $\mathcal{H}_0(W)$ is much bigger than one-dimensional subspace of singular homology classes. In fact, as was stated in our theorem, its dimension is uncountable. 

We will start with some sequence of positive numbers $n_k$ which is monotonically decreasing with $\lim n_k = 0$. From now on, up to the end of this proof, let $x_\bullet$ have a special form: 
\begin{displaymath}
x_k = (-1)^k (n_{k+1} - n_k). 
\end{displaymath}
We can see that
\begin{equation}
\sum_{k = 0}^N |x_k| = n_0 - n_{N+1}, \nonumber
\end{equation}
hence $x_\bullet \in \ell^1$. 

Let us calculate $m_k^x$ using (\ref{eq:solution}):
\begin{equation}
m_k^x = \sum_{i = 0}^k (-1)^{i+k} x_i = (-1)^k \sum_{i = 0}^k (n_{i+1} - n_i) = (-1)^k (n_{k+1} - n_0).
\label{eq:mspecialform}
\end{equation}
Since $|m_k^x|$ does not fulfil the necessary condition it is not absolutely summable. Hence, $x_\bullet$ does not correspond to the zero homology class 

More generally, we will check what conditions should $x$ satisfy in order to be a non-singular homology class. So let $y_\bullet = (\alpha, 0, 0, 0, ...)$ (for $\alpha \in \mathbb{R}$) be a sequence corresponding to some singular homology class. In this case obviously:
\begin{equation}
m_k^{x-y} = (-1)^k (n_{k+1} - n_0 - \alpha); \nonumber
\end{equation}
we can easily see this when we notice that $m^x_\bullet$ is linear with respect to $x$ according to equation (\ref{eq:solution}). So, if we take $\alpha = -n_0$ the sequence satisfies the necessary condition of series convergence. Then, we see that a sufficient condition for $x$ to be a non-singular homology class is
\begin{equation}
\sum_{k = 0}^\infty n_k = \infty, \nonumber
\end{equation}
so we are interested in sequences converging to zero but not too fast. 

As an example of such sequence we consider:
\begin{displaymath}
n_k^\beta = \frac{1}{(k+1)^\beta},
\end{displaymath}
with $0 < \beta < 1$. Now, we shall prove that the homology classes in $\mathcal{H}_0(W)$ corresponding to sequences $x^\beta_\bullet$ such that $x_k^{\beta} = (-1)^k (n_{k+1}^{\beta} - n_k^{\beta})$ form a set of linearly independent vectors. So take a finite sequence of numbers $0 < \beta_i < 1$ in an increasing order, and some finite sequence of real numbers $b_i$. We shall prove that the homology class of $z_\bullet = \sum_i b_i x_\bullet^{\beta_i}$ is nontrivial.

In order to do this we need to prove that the sequence
\begin{displaymath}
m_k^z =(-1)^k  \sum_i b_i \left( \frac{1}{(k+2)^{\beta_i}} - 1 \right)
\end{displaymath}
is not absolutely summable. To obtain the above formula we use the fact that $m_\bullet^x$ is linear with respect to $x$, and the equation (\ref{eq:mspecialform}). 

First, we notice that for the necessary condition of convergence for series $\sum_{k=0}^\infty |m_k^z|$ to be satisfied, we should have $\sum_i b_i = 0$. Then, the study of the absolute summability of the above sequence can be reduced to the study of 
\begin{displaymath}
\sum_{k=0}^\infty \left| \sum_i \frac{ b_i}{(k+2)^{\beta_i}}  \right|.
\end{displaymath}
For sufficiently big $k$ the expression in $|\cdot |$ has the sign of $b_0$ (since $\beta_0$ is the smallest of the numbers), so we can consider: 
\begin{displaymath}
\sum_{k=0}^\infty \sum_i \frac{b_i}{(k+2)^{\beta_i}}.
\end{displaymath}
This series is divergent. The easiest way to see this is to use integral criterion. First, we need to notice, that it is for monotonic sufficiently big $k$. Then, the application of the criterion is straightforward.
\begin{flushright} $\square$ \end{flushright}

\section[Berlanga topology]{Berlanga topology on Milnor-Thurston homology groups}
\label{sec:berlanga1}

Berlanga considered Milnor-Thurston homology groups equipped with some natural topology \cite{Berla}. He showed that these homology groups are functors from the category of second countable separable topological spaces to the category of locally convex topological vector spaces (not necessarily Hausdorff!).

Let $X$ be a second countable separable topological space and let $k$ be a non-negative integer. Given any function $f: C^0(\Delta^k, X) \rightarrow \mathbb{R} $, we start with a linear functional defined for $\mu \in \mathcal{C}_k(X)$ in the following way
\begin{displaymath}
\Lambda_f(\mu) = \int_{C^0(\Delta^k, X)} f d\mu.
\end{displaymath} 
We shall work with the weakest topology on $\mathcal{C}_k(X)$ such that all functionals of this type are continuous. 

Berlanga proved that if $X$ is a second countable separable topological space then each $\mathcal{C}_k(X)$ equipped with this weak topology is a locally convex Hausdorff topological vector space \cite[Assertion 2.2]{Berla} and every induced map (including boundary operator) is continuous \cite[Assertion 2.1]{Berla}. Thus, $\mathcal{C}_*$ is a functor from the category of second countable separable topological spaces to the category of locally convex Hausdorff topological vector spaces. Moreover, cycles $\mathcal{Z}_k(X)$ and boundaries $\mathcal{B}_k(X)$ with induced topology also satisfy these conditions. Therefore, homology groups 
\begin{displaymath}
\mathcal{H}_k(X) = \mathcal{Z}_k(X) / \mathcal{B}_k(X)
\end{displaymath}
can be endowed with the structure of locally convex topological vector space as quotients of locally convex topological vector spaces. This structure shall be called \emph{Berlanga topology} on $\mathcal{H}_k(X)$.

The question posed by Berlanga in \cite{Berla} is whether Milnor-Thurston homology groups are Hausdorff in this topology. As we mentioned above cycles $\mathcal{Z}_k(X)$  and boundaries $\mathcal{B}_k(X)$ are both Hausdorff, hence their quotient will be Hausdorff if and only if $\mathcal{B}_k(X)$ is closed in $\mathcal{Z}_k(X)$. 

There are two results in this direction. Firstly, Berlanga's paper that was mentioned above ends with a proof that $\mathcal{H}_1$ is always Hausdorff for spaces that are homotopy equivalent to countable CW-complexes. Secondly, Zastrow constructed an example of a space $V$ where $\mathcal{H}_0(V)$ is not Hausdorff \cite{ZastrConf}. This space $V$ is the Warsaw Circle with a part of the accumulation line removed. Based on different arguments than in \cite{ZastrConf} we obtain result for the Warsaw Circle itself: 

\begin{thm}
The Milnor-Thurston homology group $\mathcal{H}_0(W)$ is not Hausdorff in Berlanga topology.
\end{thm}

\textbf{Proof.} In order to prove that $\mathcal{H}_0(W)$ is not Hausdorff we shall find a sequence of boundaries such that the limit of this sequence is not a boundary. From the proof of Theorem \ref{thm:undim} we know that chains in $\mathcal{C}_0(W)$ can be described by elements of $\ell^1$. So, let our sequence of chains be described by elements $x_\bullet^n \in \ell^1$ in the following way:
\begin{equation}
x_k^n = \left\{ \begin{array}{lr}
- \sum_{i=1}^n (-1)^i (n_{i+1} - n_i), & \textrm{for } k = 0, \\
(-1)^k(n_{k+1} - n_k), & \textrm{for } 0 < k \leq n, \\
0, & \textrm{for } k > n. 
\end{array} \right.\nonumber
\end{equation}
where $n_\bullet \notin  \ell^1$ is a decreasing sequence of positive numbers converging to zero (compare with proof of Theorem \ref{thm:undim}).

Chains described by $x_\bullet^n$ are boundaries (or, equivalently, they represent zero homology classes). To justify it, recall the proof of Theorem \ref{thm:main1}.  We saw that an arbitrary sequence $z_\bullet \in \ell^1$ with at most $N$ nonzero elements represents the same homology class as the sequence $(\alpha, 0, 0, \dots)$, where $\alpha = \sum_{k = 0}^N (-1)^k z_k$. Therefore, we see that for each $n$ the sequence $x_\bullet^n$ represents the zero homology class. 

The natural candidate for the limit of $x_\bullet^n$ is:
\begin{equation}
x_k = \left\{ \begin{array}{lr}
- \sum_{i=1}^\infty (-1)^i (n_{i+1} - n_i), & \textrm{for } k = 0, \\
(-1)^k(n_{k+1} - n_k), & \textrm{for } k > 0.
\end{array} \right. \nonumber
\end{equation}
In order to show that the above sequence is the limit of $x_\bullet^n$ we need to prove that
\begin{displaymath}
\lim_{n \rightarrow \infty} \int_W f d(\mu - \mu_{n}) = 0, 
\end{displaymath}
for any continuous $f: W \rightarrow \mathbb{R}$. Here $\mu$ and $\mu_{n}$ are measures on $W$ corresponding to $x_\bullet$ and $x_\bullet^n$ respectively (remember that we identify $C^0(\Delta^0, W)$ with $W$). 

The measures $\mu$ and $\mu_{n}$ are concentrated on a countable set of points (it is maxima $\overline{u}_k$ and minima $\overline{l}_k$ of the sinusoid), therefore the above integral can be calculated as an infinite series. Values of the continuous function $f$ on that countable set of points form a bounded sequence $a_\bullet$, so we need to prove that
\begin{displaymath}
\lim_{n \rightarrow \infty} \left(- a_0 \sum_{i = n+1}^\infty (-1)^i (n_{i+1} - n_i) + \sum_{i = n+1}^\infty (-1)^i a_i (n_{i+1} - n_i)  \right) = 0.   
\end{displaymath}
We can easily see that it is true since tails of absolutely convergent series converge to zero. 

Assume that the homology class described by $x_\bullet$ is a boundary. Let $y_k = (-1)^k(n_{k+1} - n_k)$. Then, consider the difference
\begin{equation}
y_k - x_k = \left\{ \begin{array}{lr}
\sum_{i=0}^\infty (-1)^i (n_{i+1} - n_i), & \textrm{for } k = 0, \\
0, & \textrm{for } k > 0.
\end{array} \right.\nonumber
\end{equation}
We assumed that on the level of homology $x_\bullet$ represents zero, and thus it represents a singular homology class. On the other hand, from the above equation we see that $y_\bullet - x_\bullet$ also represents a singular homology class. Therefore, $y_\bullet$ should also represent a singular homology class. However, $y_\bullet$ is exactly the form of a sequence considered in the proof of Theorem \ref{thm:undim}, and we know that it represents a non-singular homology class (note that the sequence denoted here by $y_\bullet$ was denoted by $x_\bullet$ in the proof of that theorem). Hence, we got a contradiction. Consequently, we see that $x_\bullet$ is not a boundary and $\mathcal{H}_0(W)$ is not Hausdorff. 
\begin{flushright}$\square$\end{flushright}






\begin{thebibliography}{00}

\bibitem{SierpGask} S. Akiyama, G. Dorfer, J.M. Thuswaldner, R. Winkler, \emph{On the fundamental group of the Sierpinski-gasket}, Topology Appl. 156 (2009), 1655 -- 1672.

\bibitem{Berla} R. Berlanga, \emph{A topologised measure homology}, Glasgow Math. J. 50(2008) 359 -- 368

\bibitem{BogSier} W.A. Bogley, A.J. Sieradski, \emph{Universal path spaces}, preprint availble at: http://oregonstate.edu/ bogleyw/

\bibitem{Grif} O. Bogopolski, A. Zastrow, \emph{The word problem for some uncountable groups given by countable words}, Topology Appl. 159 (2012) 569 -- 586

\bibitem{CanCon} J. W. Cannon,  G. R. Conner, \emph{The combinatorial structure of the Hawaiian Earring group}, Topology Appl. 106 (2000) 225 -- 271

\bibitem{KEdaTalk} K. Eda, \emph{Group theoretic properties for wild algebraic topology}, talk at the Worshop on Topology of Wild Spaces and Fractals (July 4 - 8), 2011, Strobl, Austria, abstract available at http://dmg.tuwien.ac.at/dorfer/wild\_topology/abstracts.pdf

\bibitem{KEda} K. Eda, K. Kawamura, \emph{The singular homology of the Hawaiian Earring}, J. London Math. Soc. 62 (2000) 305 -- 310

\bibitem{SKHansen} S.K. Hansen, \emph{Measure homology}, Math. Scand. 83(1998) 205--219 

\bibitem{DRepovs} U. H. Karimov, D. Repov\v{s}, \emph{On the homology of the Harmonic Archipelago} Cent. Eur. J. Math. 10:3 (2012) 863 -- 872

\bibitem{Claraloh} C. L\"{o}h, \emph{Measure homology and singular homology are isometrically isomorphic}, Math. Z. 253(2006) 197 -- 218 

\bibitem {MorgMor} J. W. Morgan, I. A. Morrison, \emph{A Van Kampen theorem for weak joins}, Proc. London Math. Soc. 53 (1986) 562 -- 576.

\bibitem{Ratcl} J.G. Ratcliffe, \emph{Foundations of Hyperbolic Manifolds}, Graduate Texts in Mathematics No. 149, Springer-Verlag, New York-Heidelberg-Berlin, 1994.

\bibitem{Thurst} W.P. Thurston, \emph{Geometry and Topology of Three-manifolds}, Lecture notes, Available at http://www.msri.org/publications/books/gt3m, Princeton, 1978

\bibitem{AZNontame} A. Zastrow, \emph{On recent developments and concepts in the algebraic topology of non-tame spaces}, talk at the conference (July 6 - 11, 2009), Li\`ege, Belgium, abstract available at http://atlas-conferences.com/cgi-bin/abstract/select/cafv-01?session=1
 
\bibitem{AZastrow} A. Zastrow, \emph{On the (non)-coincidence of Milnor-Thurston homology theory with singular homology theory}, Pacific J. Math. 186(1998) 369 -- 396

\bibitem{ZastrConf} A. Zastrow, \emph{The non-Hausdorffness of Milnor-Thurston homology group}, talk at the Conference on Algebraic Topology CAT'09 (July 6 - 11, 2009), Warsaw, Poland, abstract available at http://www.mimuw.edu.pl/$\sim$cat09/abstracts.pdf



\end{thebibliography}



\end{document}